\documentclass[11pt,english]{article}
\usepackage[T1]{fontenc}

\usepackage{graphicx, subfigure}
\expandafter\def\csname ver@subfig.sty\endcsname{}
\usepackage[latin9]{inputenc}
\usepackage{geometry}
\usepackage{amsmath}
\usepackage{appendix}
\usepackage{amssymb}
\usepackage{esint}
\usepackage{makecell,multirow,diagbox}
\usepackage{mathtools}
\usepackage{epstopdf}
\usepackage[latin9]{inputenc}
\usepackage{pict2e}
\usepackage{xcolor}
\usepackage{colortbl,booktabs}
\usepackage{stmaryrd}
\usepackage{esint}
\usepackage[all,cmtip]{xy}
\usepackage{mathtools}
\usepackage{float}
\usepackage{setspace}
\usepackage{authblk}
\usepackage{amsthm}
\usepackage{mathrsfs}
\usepackage{stmaryrd}
\usepackage{bbold}
\usepackage{color}
\usepackage{tikz}
\usepackage[all,cmtip]{xy}
\usepackage{algorithm}  
\usepackage{algorithmic}

\theoremstyle{plain}\newtheorem{theorem}{Theorem}[section]
\theoremstyle{plain}
\newtheorem{remark}{Remark}
\usepackage{color}
\definecolor{marin}{rgb} {0., 0.3, 0.7}
\definecolor{rouge}{rgb} {0.8, 0., 0.}
\definecolor{sepia}{rgb} {0.8, 0.5, 0.}
\usepackage[colorlinks,citecolor=sepia,linkcolor=marin,bookmarksopen,bookmarksnumbered]{hyperref}

\def\be{\begin{equation}}
\def\ee{\end{equation}}

\theoremstyle{definition}

\DeclareSymbolFont{largesymbol}{OMX}{yhex}{m}{n}
\DeclareMathAccent{\Widehat}{\mathord}{largesymbol}{"62}

\makeatletter

\usepackage{babel}
\usepackage{subfig}
\graphicspath{{Dirac/}}
           
\makeatother
\DeclareMathSizes{14}{14}{9.8}{7}

\usepackage{lipsum}

\begin{document}

\title{Symplectic particle-in-cell methods for hybrid plasma models with Boltzmann electrons and space-charge effects}
\date{}
\author{Yingzhe Li\thanks{yingzhe.li@ipp.mpg.de}}
\affil{Max Planck Institute for Plasma Physics, Boltzmannstrasse 2, 85748 Garching, Germany}

\maketitle
\begin{abstract}
We study the geometric particle-in-cell methods for an electrostatic hybrid plasma model. In this model, ions are described by the fully kinetic equations, electron density is determined by the Boltzmann relation, and space-charge effects are incorporated through the Poisson equation. By discretizing the action integral or the Poisson bracket of the hybrid model, we obtain a finite dimensional Hamiltonian system, for which the Hamiltonian splitting methods or the discrete gradient methods can be used to preserve the geometric structure or energy.
The global neutrality condition is conserved under suitable boundary conditions. Moreover, the results are further developed for an electromagnetic hybrid model proposed in [Vu H X. J Comput Phys, 124(2):417-430]. Numerical experiments of finite grid instability, Landau damping, and resonantly excited nonlinear ion waves illustrate the behaviour of the proposed numerical methods.
\end{abstract}
\setcounter{tocdepth}{1} 


\section{Introduction}        
Hybrid plasma models with Boltzmann electrons and space-charge effects~\cite{vu1996adiabatic, cartwright2000nonlinear, tajima2018computational} constitute an important class of plasma models. In these models, the electron density is directly determined from the potential via the Boltzmann relation, and space-charge effects are included via the Poisson equation. 
The electrostatic hybrid model with Boltzmann electrons and space-charge effects (HBS model) has many applications in plasma physics. The acceleration of light and heavy ions in an expanding plasma slab with hot electrons produced by an intense and short laser pulse is studied using the HBS model in~\cite{ionacc}. In~\cite{huwang}, by numerical simulations of the HBS model, the expansion of a collisionless hypersonic plasma plume into a vacuum is investigated.
In~\cite{cohen1997resonantly} resonantly excited nonlinear ion waves are investigated numerically using the HBS model, and it is noted that the exponential term in the Poisson equation introduces sufficient nonlinearity, allowing us  to derive the dispersion relation for parametric instabilities and describe the generation of the second harmonic. Mathematically, a related model with Boltzmann electrons is derived and proved to be well-posed globally in~\cite{Bardos}. To include electromagnetic effects, an electromagnetic hybrid model with the self-consistent ponderomotive driving potential is proposed in~\cite{vu1996adiabatic}, and a more general fully kinetic, reduced-description particle-in-cell model is proposed in~\cite{vu1999aspen} for the ion-driven parametric instabilities.

Different from the fully kinetic (for both ions and electrons) Vlasov--Poisson system, there is no electron distribution function in the HBS model. The simulations using the HBS model allow time step sizes on the scale of ions and are thus more efficient.
By taking the quasi-neutral limit of the HBS model, a simplified hybrid model without the space-charge effects~\cite{rambo1995finite} can be derived. However, the space-charge effects are important and needed to be incorporated in some cases~\cite{vu1996adiabatic}, such as in the inertial confined fusion regime where $k\lambda_e = \mathcal{O}(1)$, with $\lambda_e$ the electron Debye length and $k$ the wave number. To achieve accurate resolution at the electron Debye scale, for example, to recover numerically the $k^2 \lambda_e^2$ term in the dispersion relation of the ion acoustic waves, the mesh size $\Delta x$ must satisfy $\Delta x < \lambda_e$.

There have been a lot of numerical methods developed for electrostatic plasma models, such as Eulerian methods~\cite{mdv, discontinuous}, particle-in-cell methods~\cite{pic1, pic2}, and semi-Lagrangian methods~\cite{Cheng, semi}. Recently, some structure-preserving methods have been proposed in~\cite{webb, sunvp} for the fully kinetic Vlasov--Poisson system. Structure-preserving methods for the hybrid model with quasi-neutrality and Boltzmann electrons have been developed based on variational or Hamiltonian formulations~\cite{drift, fB, fA}. 
The nonlinear Poisson--Boltzmann equation in the HBS model appears in many electrostatic models in biomolecular simulations. For a review about fast analytical methods, see~\cite{zhenli}, and for numerical methods, see~\cite{lu}.  Plenty of numerical methods have been proposed for the Poisson--Boltzmann equation, such as the finite element method~\cite{chenlong} and the iterative discontinuous Galerkin method~\cite{huang1, huang2}. 

 In this work, our discretizations of the HBS model follow the structure-preserving methods for models in plasma physics~\cite{DECVM, Qincanonical, he2016hamiltonian, GEMPIC, morrison2017structure}, which preserve the geometric structures of the systems and exhibit very good long-term behavior~\cite{Feng, HLW}. 
The numerical schemes constructed in this work complement existing structure-preserving methods for other (hybrid) electrostatic models. Moreover, a more complicated electromagnetic hybrid model proposed in~\cite{vu1996adiabatic} is investigated. 

The action integral and Hamiltonian structure of the HBS model in this work are derived based on the results in~\cite{low1958lagrangian, drift, morrison1980maxwell}. 
By discretizing the action integral as in~\cite{DECVM, drift} or the Poisson bracket as in~\cite{Qincanonical,GEMPIC} with particle methods for the distribution function and finite difference methods for the electrostatic potential, we obtain a finite dimensional Hamiltonian system. Time discretizations are conducted using the Hamiltonian splitting methods~\cite{HLW} and the discrete gradient methods~\cite{mclachlan1999geometric, Gonzalez}. 
In plasma physics simulations, Hamiltonian splitting methods have been used in~\cite{crouseilles2015hamiltonian, qin2015comment, he2015hamiltonian} and discrete gradient methods have been used in~\cite{kormann2021energy} as time integrators. For the electromagnetic hybrid model~\cite{vu1996adiabatic}, a Poisson bracket is proposed as the sum of the Lie--Poisson bracket~\cite{morrison1980maxwell} and the canonical bracket of the Schr\"odinger equation~\cite{marsden2013introduction}.

 The neutrality condition is preserved by the discretizations of the HBS model with appropriate boundary conditions. Moreover, we demonstrate that the quasi-neutral limits of the schemes proposed are structure-preserving for the hybrid model with quasi-neutrality and Boltzmann electrons. The numerical methods are validated by the good conservation of energy. We conduct the implementation in the Python package STRUPHY~\cite{STRUPHY}. 
 
The paper is organized as follows. In Section~\ref{sec:model}, the action integral and the Poisson bracket are presented for the HBS model.  In Section~\ref{sec:dis}, structure preserving discretizations are given. In Section~\ref{sec:neutrallimit}, two asymptotic limits and the dispersion relation of the linear Landau damping of the HBS model are discussed. Geometric structure and discretization of the electromagnetic hybrid model proposed in~\cite{vu1996adiabatic} are presented in Section~\ref{sec:extension}. In Section~\ref{sec:numsec}, numerical experiments of finite grid instability, Landau damping, and resonantly excited nonlinear ion waves are conducted to validate the numerical schemes of the HBS model.  In Section~\ref{sec:conclusion}, we conclude the paper with a summary and an outlook for future works.

 \section{The electrostatic hybrid plasma model with Boltzmann electrons and space-charge effects}
\label{sec:model}
In this section, we introduce the action integral and the Poisson bracket for the HBS model, and formulate the model as a Hamiltonian system. The electromagnetic hybrid model proposed in~\cite{vu1996adiabatic} is presented in section~\ref{sec:extension}.
The HBS model with physical units is 
\begin{equation*}
\begin{aligned}
&\frac{\partial f}{\partial t} + {\mathbf v} \cdot \nabla f +\frac{Ze}{m_i}{\mathbf E} \cdot \nabla_v f = 0,\\
& {\mathbf E} = - \nabla \phi,\\
&-\epsilon_0 \Delta \phi = Z e \int f\, \mathrm{d}{\mathbf v} - e n_0\exp\left({\frac{e(\phi-\phi_0)}{k_B T_e}}\right), \quad \text{Poisson--Boltzmann}.
\end{aligned}
\end{equation*} 
Here $f(t, {\mathbf x}, {\mathbf v})$ represents the distribution function of ions, which depends on time $t$, position ${\mathbf x}$, and velocity ${\mathbf v}$. The symbol $e$ denotes the unit charge, $m_i$ denotes the ion mass, $Z$ denotes the ion charge number. The electrostatic potential 
$\phi(t, {\mathbf x})$ is  determined by the Poisson--Boltzmann equation, and electron density $n_e$ is related to $\phi$ through the Boltzmann relation,
$$n_e = n_0 \exp({(\phi-\phi_0)/T_e}), \quad \text{Boltzmann relation},$$
where $n_0(t,{\mathbf x})$ is the reference electron number density, $\phi_0(t, {\mathbf x})$ is the reference potential or the low-frequency ponderomotive potential, and $T_e(t, {\mathbf x})$ is the given temperature of electrons. 

The normalization is done as
\begin{equation}\label{eq:llll}
\tilde{\mathbf x} = \frac{\mathbf x}{\lambda_D}, \ \tilde{\mathbf v} = \frac{\mathbf v}{C_0}, \ \tilde{t} = t \omega_i, \ \tilde{f} = \frac{C_0^3}{\bar{n}}f, \ \tilde{n}_0 = \frac{n_0}{\bar{n}}, \ \tilde{T}_e = \frac{T_e}{\bar{T}_i}, \ \tilde{\phi} = \frac{e\phi}{k_B \bar{T}_i}, \ \tilde{\phi}_0 = \frac{e\phi_0}{k_B \bar{T}_i},
\end{equation} 
where $C_0 = \sqrt{\frac{k_B \bar{T}_i}{m_i}}$ is the ion thermal speed, $Z\omega_i = \sqrt{\frac{\bar{n} Z^2 e^2}{\epsilon_0 m_i}}$ is the ion plasma frequency, $\lambda_D = C_0/\omega_i = \sqrt{\frac{\epsilon_0 k_B \bar{T}_i}{\bar{n} e^2}}$, $\bar{n}$ is the characteristic ion density, and $\bar{T}_i$ is the characteristic ion temperature. Then we get the normalized HBS model (tilde symbol is omited for convenience)
\begin{equation}\label{eq:avp}
\begin{aligned}
&\frac{\partial f}{\partial t} + {\mathbf v} \cdot \nabla f + Z{\mathbf E} \cdot \nabla_v f = 0,\\
& {\mathbf E} = - \nabla \phi,\\
&-\Delta \phi = Z \int f\, \mathrm{d}{\mathbf v} - n_0 \exp({(\phi-\phi_0)/T_e}), \quad \text{Poisson--Boltzmann}.
\end{aligned}
\end{equation} 
When periodic or zero Neumann boundary condition is imposed, The HBS model satisfies a neutrality condition given by,
 $$Z\int f\, \mathrm{d}{\mathbf x}\mathrm{d}{\mathbf v} = \int n_0 \exp({(\phi-\phi_0)/T_e}) \mathrm{d}{\mathbf x}.$$
 To construct structure-preserving methods for this model, the following action integral and Poisson bracket are proposed. For convenience, we consider the case with time independent  $n_0, \phi_0, T_e$. The time dependent case can be addressed using the technique of extending the dimension~\cite{zhou2017explicit}. \\
\noindent{\bf Variational principle} By adding the term $\frac{1}{2} |\nabla \phi|^2$ in the Low's action principle~\cite{low1958lagrangian}, and combining it with the action principle proposed in~\cite{drift}, we derive the following action integral
\begin{equation}\label{eq:vpai}
\mathcal{A}({\mathbf x}, \phi) = \int f_0({\mathbf x}_0, {\mathbf v}_0) \left( \frac{|\dot{\mathbf x}|^2}{2} - Z\phi({\mathbf x})\right) \mathrm{d}{\mathbf z}_0 + \int \frac{|\nabla \phi|^2}{2} \mathrm{d}{{\mathbf x}}+ \int n_0 T_e \exp\left(\frac{\phi-\phi_0}{T_e}\right)\mathrm{d}{{\mathbf x}},
\end{equation}
where $\mathrm{d}{{\mathbf z}_0} := \mathrm{d}{{\mathbf x}_0}\mathrm{d}{{\mathbf v}_0}$, ${\mathbf x} ={\mathbf x}({\mathbf x}_0, {\mathbf v}_0, t)$, and $\dot{\mathbf x} =\mathrm{d}{\mathbf x}({\mathbf x}_0, {\mathbf v}_0, t)/\mathrm{d}t$. We introduce ${\mathbf v} = \dot{\mathbf x} $ and  $f(t, {\mathbf x}, {\mathbf v}) = f_0({\mathbf x}_0, {\mathbf v}_0)$, 
the Euler--Lagrangian equations 
$\frac{\delta \mathcal{A}}{\delta {\mathbf x}} = 0, \frac{\delta \mathcal{A}}{\delta \phi} = 0$
can be written as 
$$
\ddot{\mathbf x} = -Z\nabla \phi(\mathbf x), \quad -\Delta \phi = Z \int f(t, {\mathbf x}, {\mathbf v}) \mathrm{d}{\mathbf v} - n_0 \exp({(\phi-\phi_0)/T_e}),
$$
which yields the HBS model~\eqref{eq:avp} by calculating ${\mathrm{d}f}/{\mathrm{d}t} = 0$.

\noindent{\bf Poisson bracket} The Poisson bracket of this model is the same as the Vlasov--Poisson system's Poisson bracket proposed in~\cite{morrison1980maxwell},
\begin{equation}\label{eq:liepoisson}
\{ \mathcal{F}, \mathcal{G} \}(f) = \int f \left[\frac{\delta \mathcal{F}}{\delta f},  \frac{\delta \mathcal{G}}{\delta f}\right]_{xv} \mathrm{d}{\mathbf x}\mathrm{d}{\mathbf v},
\end{equation}
where $[g, h]_{xv} = \nabla_{\mathbf x} g \cdot \nabla_{\mathbf v} h - \nabla_{\mathbf x} h \cdot \nabla_{\mathbf v} g$.
The Hamiltonian (total energy) of this model is
\[
\begin{aligned}
\mathcal{H} & = \int |{\mathbf v}|^2 f \mathrm{d}{\mathbf x}\mathrm{d}{\mathbf v} - \mathcal{A}\\
& = - \frac{1}{2} \int |\nabla \phi|^2 \mathrm{d}{\mathbf x} + Z \int f \phi \,\mathrm{d}{\mathbf x}\mathrm{d}{\mathbf v} - \int T_e n_0 \exp\left(\frac{\phi-\phi_0}{T_e}\right) \mathrm{d}{\mathbf x} + \frac{1}{2} \int |{\mathbf v}|^2 f\, \mathrm{d}{\mathbf x}\mathrm{d}{\mathbf v}.
\end{aligned}
\]
Based on the bracket and Hamiltonian, the HBS model~\eqref{eq:avp} can be formulated as 
$$
\dot{f} = \{ f, \mathcal{H} \}.
$$ 
Here we regard $f$ as the only unknown of the HBS model~\eqref{eq:avp}, and $\phi$ is determined by $f$ from the Poisson--Boltzmann equation in~\eqref{eq:avp}.

\section{Discretization}\label{sec:dis}
We use the particle-in-cell methods to discretize the distribution function $f$ and finite difference methods to discretize the electrostatic potential $\phi$. Two equivalent structure-preserving phase-space discretizations are obtained by discretizing the action integral and the Poisson bracket. 
The Hamiltonian splitting method and the discrete gradient method are used for time discretizations to preserve the geometric structure and energy, respectively. In the following, $a^n$ and $a^{n+1}$ represent the  approximations of $a$ at times $n\Delta t$ and $(n+1)\Delta t$, respectively, and 
$a^{n+\frac{1}{2}} = \frac{ a^{n} + a^{n+1}  }{2}$, where $\Delta t$ is the time step size.

\subsection{Discretization of $f$ and $\phi$}
Here we focus on the one dimensional case with periodic boundary condition, higher dimensional cases can be treated similarly.
The distribution function $f$ is approximated as
$$
 f_h({x}, {v}, t) = \sum_{k=1}^{N_p} w_k S ({x} - {x}_k) \delta ({v} - {v}_k),
$$
where $N_p$ is the total particle number, $w_k$, ${x}_k$, and ${v}_k$ denote the weight, position, and velocity of $k$-th particle. $S$ is the shape function of particle,  typically chosen as a B-spline. We use the vector ${\mathbf X}$ to denote $(x_1, \cdots, x_{N_p})^\top$, and vector ${\mathbf V}$ to  denote $(v_1, \cdots, v_{N_p})^\top$.

The electrostatic potential $\phi$ is discretized using finite difference method, i.e. 
$$
\phi_j \approx\phi({x}_j), \quad j = 1, \cdots, N,
$$
with a set of uniform grids $ \{ {x}_j \} $, $N$ is the number of grids, 
$
(\phi_1, \cdots, \phi_N)^\top
$
is denoted as ${\boldsymbol \phi}$, and $(\phi_0(x_1), \cdots, \phi_0(x_N))^\top $ is denoted as ${\boldsymbol \phi}_0$.

\subsection{Phase-space discretization}
\subsubsection{Discretization of action integral}
We approximate variational action integral~\eqref{eq:vpai} as 
\begin{equation*}
\begin{aligned}
\mathcal{A}_h({\mathbf X}, {\boldsymbol \phi}) & = \sum_{k=1}^{N_p} w_k \left( \frac{1}{2} \dot{x}_k^2 - Z \sum_{j=1}^N \Delta x  S({x}_j - {x}_k) \phi_j \right) + \frac{1}{2}  {\boldsymbol \phi}^\top \mathbb{A}   {\boldsymbol \phi} \Delta x\\
& + \sum_{j=1}^N \Delta x n_0({x}_j) T_e({x}_j) \exp\left(\frac{\phi_j - \phi_0(x_j)}{T_e(x_j)}\right),
\end{aligned}
\end{equation*}
where the matrix $\mathbb{A}$ of size $N \times N$ is 
$$
\mathbb{A} = \frac{1}{\Delta x^2}
\left(\begin{matrix}
 -2 & 1 & 0 & \cdots & 0 & 0 & 1 \\
  1 & -2 & 1 & 0 & \cdots & 0 & 0  \\
  0 & 1 & -2 & 1 & 0 & \cdots & 0 \\
   \vdots & \ddots & \ddots  &  \ddots & \ddots & \ddots & \vdots \\
    0 & \cdots & 0 & 1 & -2 & 1 & 0\\
       0 & \cdots & 0 & 0 & 1 & -2 & 1 \\
  1 & 0 & 0 & \cdots & 0 & 1 & -2 
\end{matrix} \right). 
$$
By calculating the variations about ${x}_k$ and ${\boldsymbol \phi}$, we have 
\begin{equation}
\label{eq:disva}
\begin{aligned}
& \ddot{x}_k =  Z \sum_j \Delta x  \partial_x S({x}_j - {x}_k) \phi_j, \quad k = 1, \cdots, N_p,\\
& -Z \sum_{k=1}^{N_p} w_k  S({x}_j - {x}_k) + (\mathbb{A}   {\boldsymbol \phi})_j + n_0({x}_j) \exp\left(\frac{\phi_j - \phi_0(x_j)}{T_e({x}_j)}\right)  = 0, \quad j = 1, \cdots, N.
\end{aligned}
\end{equation}
\begin{remark}
The fixed point iteration method in~\cite{cohen1997resonantly} is used for solving the above discretized Poisson--Boltzmann equation, i.e., the second equation of~\eqref{eq:disva}, 
$$
 \mathbb{A}   {\boldsymbol \phi}^{m+1} - c{\boldsymbol \phi}^{m+1} = -{\mathbf n}_i - c{\boldsymbol \phi}^{m} + {\mathbf n}_{e}^m,
$$
where $m$ is the iteration index, the linear system in iteration $m \rightarrow m+1$ is  solved by conjugate gradient method, ${\mathbf n}_i$ is the ion density at grids, electron density at $j$-grid is $n_{e,j}^m = n_0({x}_j) \exp\left({\frac{\phi_j^m}{T_e({x}_j)}}\right)$ and $c = \text{max} \{\frac{1}{T_e(x_j)} \exp({{{\phi_j^m/T_e({x}_j), \, j = 1, \cdots, N}\}}})$. The initial value of the fixed point iteration is set as zero in Section~\ref{sec:numsec}.
\end{remark}

The equations~\eqref{eq:disva} can be formulated as a Hamiltonian system by the Legendre transformation~\cite{marsden2013introduction}. In the following, we present another way to derive a finite dimensional Hamiltonian system.

\subsubsection{Discretization of Poisson bracket}
Here we discretize the Poisson bracket~\eqref{eq:liepoisson} according to~\cite{Qincanonical, GEMPIC} and obtain 
\begin{equation}\label{eq:dispoissonbra}
\{F, G\}_h = \sum_{k=1}^{N_p} \frac{1}{w_k} \left(\partial_{x_k}F    \partial_{v_k}G -  \partial_{x_k}G    \partial_{v_k}F \right).
\end{equation}
The discrete Hamiltonian is 
\begin{equation}
\label{eq:dishamiltonian}
\begin{aligned}
H & = - \frac{1}{2}  {\boldsymbol \phi}\top \mathbb{A}   {\boldsymbol \phi} \Delta x + Z \sum_{j=1}^N \Delta x  \sum_{k=1}^{N_p} w_k S({x}_j - {x}_k) \phi_j  + \frac{1}{2} \sum_{k=1}^{N_p} w_k {v}_k^2\\
& - \sum_{j=1}^N \Delta x T_e({x}_j) n_0({x}_j) \exp\left({\frac{\phi_j - \phi_0(x_j)}{T_e({x}_j)}}\right),
\end{aligned}
\end{equation}
where the ${\boldsymbol \phi}$ is determined by the particles via the second equation in~\eqref{eq:disva}.
Then we can obtain the following finite dimensional Hamiltonian system after phase-space discretization,
\begin{equation}\label{eq:spdeha}
\dot{x}_k = \frac{1}{w_k} \partial_{v_k}H, \quad \dot{v}_k = - \frac{1}{w_k} \partial_{x_k}H, \quad k = 1, \cdots, N_p.
\end{equation} 
The following theorem shows that discretizations of the action principle and the Poisson bracket as above are equivalent. 
\begin{theorem}
The Hamiltonian system~\eqref{eq:spdeha} is equivalent to~\eqref{eq:disva}. 
\end{theorem}
\begin{proof}
As we know that $\partial_{v_k}H = w_k {v}_k$, we have $\dot{x}_k = v_k$. To obtain~\eqref{eq:disva}, the only thing we need to prove is that $\partial_{x_k}H = - Z w_k \sum_{j=1}^N \Delta x  \partial_x S({x}_j - {x}_k) \phi_j $, which is obtained by calculating $\partial_{x_k}H$ with the
discrete Poisson--Boltzmann equation (the second equation in~\eqref{eq:disva}), 
\begin{equation*}
\begin{aligned}
& \partial_{x_k}H = - Z w_k \sum_{j=1}^N \Delta x \partial_x S({x}_j - {x}_k) \phi_j - \mathbb{A} {\boldsymbol \phi} \cdot \frac{\partial \boldsymbol \phi}{\partial {x_k}} \Delta x  \\
& + Z \sum_{j=1}^N \Delta x  \sum_{k'=1}^{N_p} w_k S({x}_j - {x}_k') \frac{\partial \phi_j}{\partial {x_k}}  - \sum_{j=1}^N \Delta x_j n_0({x}_j) \exp\left({\frac{\phi_j - \phi_0(x_j)}{T_e({x}_j)}}\right) \cdot \frac{\partial \phi_j}{\partial {x_k}},
\end{aligned}
\end{equation*}
\end{proof}
where the sum of the last three terms is zero because of the discrete Poisson--Boltzmann equation (the second equation in~\eqref{eq:disva}).

\begin{theorem}{Discrete neutrality is conserved by the discretizations~\eqref{eq:disva} and~\eqref{eq:spdeha}.}
 \end{theorem}
 \begin{proof}
Taking the sum over $j$ in the discrete Poisson--Boltzmann equation (the second equation in~\eqref{eq:disva}) gives 
$$
-Z \sum_{j=1}^N \sum_{k=1}^{N_p} w_k  S({x}_j - {x}_k) \Delta x + \sum_{j=1}^N n_0({x}_j) \exp\left(\frac{\phi_j - \phi_0(x_j)}{T_e({x}_j)}\right) \Delta x = 0,
$$
 where we use that $\sum_j (\mathbb{A}   {\boldsymbol \phi})_j = 0$. This proves the discrete neutrality.
 \end{proof}

\subsection{Time discretization} 
In this subsection, we introduce the time discretizations for~\eqref{eq:disva} and~\eqref{eq:spdeha}. 
The first method is the Hamiltonian splitting method~\cite{HLW}, which is explicit and symplectic for~\eqref{eq:disva} and~\eqref{eq:spdeha}.  This method was applied in~\cite{crouseilles2015hamiltonian, qin2015comment, he2015hamiltonian} for the Vlasov--Maxwell equations, and was used in~\cite{DECVM, he2016hamiltonian, GEMPIC} for the construction of the fully structure-preserving methods. The other time discretization is the implicit discrete gradient method~\cite{mclachlan1999geometric}, which preserves the energy exactly. \\
\noindent{\bf Hamiltonian splitting method} We split the Hamiltonian~\eqref{eq:dishamiltonian} as $H = H_1 + H_2$, where 
\begin{equation}
\begin{aligned}
H_1 &= \frac{1}{2} \sum_{k=1}^{N_p} w_k {v}_k^2,\\
H_2 &= - \frac{1}{2}  {\boldsymbol \phi}^\top \mathbb{A}   {\boldsymbol \phi} \Delta x + Z \sum_{j=1}^N \Delta x  \sum_{k=1}^{N_p} w_k S({x}_j - {x}_k) \phi_j \\
& - \sum_{j=1}^N \Delta x T_e({x}_j) n_0({x}_j) \exp\left({\frac{\phi_j - \phi_0(x_j)}{T_e({x}_j)}}\right),
\end{aligned}
\end{equation}
which give the following two corresponding subsystems,
\begin{equation}
\label{eq:haad}
\begin{aligned}
& \text{sub-step I}: \dot{x}_k = {v}_k, \quad \dot{v}_k = 0\\
&  \text{sub-step II}: \dot{x}_k = 0, \quad \dot{v}_k = Z \sum_{j=1}^{N} \Delta x  \partial_x S({x}_j - {x}_k) \phi_j, 
\end{aligned}
\end{equation}
where $\phi_j \ (j = 1, \cdots, N)$ are given by the discrete Poisson--Boltzmann (the second equation in~\eqref{eq:disva}).
Both sub-steps can be solved exactly. Here we present the first and second order methods by composition method~\cite{HLW},
\begin{align*}
& \text{First order Lie splitting}:   \Phi^{1}_{\Delta t}   \circ \Phi^{2}_{\Delta t},\\
& \text{Second order Strang splitting}: \Phi^{2}_{\Delta t/2}  \circ \Phi^{1}_{\Delta t} \circ \Phi^{2}_{\Delta t/2},
\end{align*}
where $ \Phi^{1}_{\Delta t} $ and $ \Phi^{1}_{\Delta t} $ are solution maps of substeps I and II, respectively. 
Higher order structure-preserving schemes can be constructed by composition methods~\cite{HLW}.

\noindent{\bf Discrete gradient method} We use the second order discrete gradient method proposed in~\cite{Gonzalez} to converse energy exactly 
\begin{equation}
\label{eq:disenen}
\begin{aligned}
& \frac{{\mathbf X}^{n+1} - {\mathbf X}^{n}}{\Delta t} =   {\mathbb{W}}^{-1} \bar{\nabla}_{\mathbf{V}}  H,\quad \frac{{\mathbf V}^{n+1} - {\mathbf V}^{n}}{\Delta t} = -   {\mathbb{W}}^{-1} \bar{\nabla}_{\mathbf{X}}  H,
\end{aligned}
\end{equation}
where 
\begin{equation*}
\begin{aligned}
& \mathbb{W} = \text{diag}(w_1, \cdots, w_{N_p}),\\
& \bar{\nabla}_{\mathbf X} H = \nabla_{\mathbf X} H \left(\frac{{\mathbf X}^n + {\mathbf X}^{n+1}  }{2}\right) + d_{c} \left( {\mathbf X}^{n+1} - {\mathbf X}^n \right),\\
& \bar{\nabla}_{\mathbf V} H = \nabla_{\mathbf V} H \left(\frac{{\mathbf V}^n + {\mathbf V}^{n+1}  }{2}\right) + d_{c} \left( {\mathbf V}^{n+1} - {\mathbf V}^n \right),\\
& d_{c} = \frac{H_d- \nabla H ({\mathbf X}^{n+\frac{1}{2}}, {\mathbf V}^{n+\frac{1}{2}}) \cdot (  ({\mathbf X}^{n+1} - {\mathbf X}^n)^\top,   ({\mathbf V}^{n+1} - {\mathbf V}^n)^\top)^\top    }{ |{\mathbf X}^{n+1} - {\mathbf X}^n|^2 +  |{\mathbf V}^{n+1} - {\mathbf V}^n|^2 },\\
& H_d = H({\mathbf H}^{n+1}, {\mathbf V}^{n+1}) - H({\mathbf X}^n, {\mathbf V}^{n}).
\end{aligned}
\end{equation*} 
This discrete gradient method is implicit, for which fixed-point iteration method is used. The degree of shape function $S$ is chosen to be at least two for the convergence of iterations.

\section{Asymptotic limits}\label{sec:neutrallimit}

In this section we discuss two asymptotic limits, quasi-neutral limit and large $T_e$ limit, with corresponding suitable normalization. 

\subsection{Quasi-neutral limit}
We do the normalization as
\begin{equation}\label{eq:quasiscaling}
\tilde{\mathbf x} = \frac{\mathbf x}{x^*}, \ \tilde{\mathbf v} = \frac{\mathbf v}{C_0}, \ \tilde{t} = t \omega_i, \ \tilde{f} = \frac{C_0^3}{\bar{n}}f, \ \tilde{n}_0 = \frac{n_0}{\bar{n}}, \ \tilde{\phi} = \frac{e\phi}{k_B \bar{T}_e}, \ \tilde{\phi}_0 = \frac{e\phi_0}{k_B \bar{T}_e}, \ \tilde{T}_e = \frac{T_e}{\bar{T}_e},
\end{equation}
where where $x^*$ is the space scale interested, $Z\omega_i = \sqrt{\frac{\bar{n} Z^2 e^2}{\epsilon_0 m_i}}$ is the ion plasma frequency, $\lambda_D = \sqrt{\frac{\epsilon_0 k_B \bar{T}_e}{\bar{n} e^2}}$ is the electron Debye length, $C_0 = \lambda_D \omega_i =  \sqrt{\frac{k_B \bar{T}_e}{m_i}}$, $\bar{n}$ is the characteristic ion density, and $\bar{T}_e$ is the characteristic electron temperature. Then we get the normalized HBS model (tilde symbol is omited for convenience)
\begin{equation}\label{eq:quasin}
\begin{aligned}
&\frac{\partial f}{\partial t} + {\mathbf v} \cdot \nabla f + Z{\mathbf E} \cdot \nabla_v f = 0,\\
& {\mathbf E} = - \nabla \phi,\\
&- \lambda^2 \Delta \phi = Z \int f\, \mathrm{d}{\mathbf v} - n_0 \exp\left(\frac{\phi-\phi_0}{T_e}\right), \quad \text{Poisson--Boltzmann},
\end{aligned}
\end{equation} 
where $\lambda = \lambda_D / x^*$.

When we take the quasi-neutral limit $\lambda \rightarrow 0$ for the normalized HBS model~\eqref{eq:quasin},
 we get the following hybrid model with quasi-neutrality and Boltzmann electrons~\cite{rambo1995finite, drift},
\begin{equation}\label{eq:pq}
\begin{aligned}
&\frac{\partial f}{\partial t} + {\mathbf v} \cdot \nabla f + Z{\mathbf E} \cdot \nabla_v f = 0,\\
& {\mathbf E} = - \nabla \phi,\\
&0 = Z \int f \,\mathrm{d}{\mathbf v} - n_0 \exp\left(\frac{\phi - \phi_0}{T_e}\right).
\end{aligned}
\end{equation} 
By the similar calculations for the Vlasov--Poisson equation in~\cite{ericvlasov}, we get 
the dispersion relation of the linear Landau damping of~\eqref{eq:quasin} with $Z = 1, n_0=1, \phi_0 = 0$
\begin{equation}
\label{eq:erichbsdisper}
1 + \lambda^2{k^2 T_e} = \frac{T_e}{T_i}  \mathcal{Z}'\left(\frac{\omega}{kv_T} \right),
\end{equation}
where $ \mathcal{Z}$ is the plasma dispersion function and $T_i = v^2_T/2$.
By $\lambda \rightarrow 0$ we get the dispersion relation of~\eqref{eq:pq}~\cite{kunz2014pegasus}
$$
1 = \frac{T_e}{T_i}   \mathcal{Z}'\left(\frac{\omega}{kv_T} \right).
$$

Under this normalization~\eqref{eq:quasiscaling}, the discrete Poisson--Boltzmann equation in our scheme~\eqref{eq:haad} becomes 
\begin{equation}
\label{eq:quasidispb}
\begin{aligned}
-Z \sum_{k=1}^{N_p} w_k  S({x}_j - {x}_k) + \lambda^2(\mathbb{A}   {\boldsymbol \phi})_j + n_0({x}_j) \exp\left(\frac{\phi_j - \phi_0(x_j)}{T_e({x}_j)}\right)  = 0, \quad j = 1, \cdots, N.
\end{aligned}
\end{equation}
When taking quasi-neutral limit for the scheme~\eqref{eq:haad} with~\eqref{eq:quasidispb},  
we get the following structure-preserving scheme similar to the scheme proposed in~\cite{drift} derived using discrete exterior calculus and Whitney form,
\begin{equation}
\label{eq:halimit}
\begin{aligned}
& \text{sub-step I}: \dot{x}_k= {v}_k, \quad \dot{v}_k= 0\\
&  \text{sub-step II}: \dot{x}_k = 0, \quad \dot{v}_k = Z \sum_{j=1}^N \Delta x  \partial_x S({x}_j - {x}_k) \phi_j, 
\end{aligned}
\end{equation}
where $\phi_j, j = 1, \cdots, N$ are determined by the following discrete equation about ${\boldsymbol \phi}$
$$
 -Z \sum_{k=1}^{N_p} w_k  S({x}_j - {x}_k) + n_0({x}_j) \exp\left(\frac{\phi_j - \phi_0(x_j)}{T_e}\right)  = 0, \quad j = 1, \cdots, N.
$$
Then limiting scheme~\eqref{eq:halimit} is the Hamiltonian splitting method for quasi-neutral limit model~\eqref{eq:pq}, i.e., scheme~\eqref{eq:haad}, is asymptotic preserving~\cite{ap} and structure-preserving at the same time. 
Similarly, the discrete gradient method~\eqref{eq:disenen} for the HBS model becomes a discrete gradient method for the hybrid model with quasi-neutrality and Boltzmann electrons.  Note that quasi-neutral limit is not a singular asymptotic limit as explained in~\cite{degond} for the case of Euler-Poisson-Boltzmann model.

\subsection{Large $T_e$ limit}\label{sec:largetelimit}
Here we adopt the normalization~\eqref{eq:llll}. By taking the $T_e \rightarrow +\infty$ in~\eqref{eq:avp}, we get the equations 
\begin{equation}\label{eq:ffvp}
\begin{aligned}
&\frac{\partial f}{\partial t} + {\mathbf v} \cdot \nabla f + Z{\mathbf E} \cdot \nabla_v f = 0,\\
& {\mathbf E} = - \nabla \phi,\\
&-\Delta \phi  = Z \int f\, \mathrm{d}{\mathbf v} - n_0.
\end{aligned}
\end{equation} 
Divided by $k^2T_e$ on the both sides of the dispersion relation~\eqref{eq:erichbsdisper} of~\eqref{eq:avp}, we get the the following dispersion relation with the current normalization
$$
1 + \frac{1}{k^2 T_e} = \frac{1}{k^2 T_i}   \mathcal{Z}'\left(\frac{\omega}{kv_T} \right).
$$
By $T_e \rightarrow +\infty$ we get the dispersion relation of model~\eqref{eq:ffvp}~\cite{ericvlasov}
$$
1 = \frac{1}{k^2 T_i}   \mathcal{Z}'\left(\frac{\omega}{kv_T} \right).
$$
Similar to the quasi-neutral limit, when $T_e \rightarrow +\infty$, the limiting schemes of the Hamiltonian splitting method~\eqref{eq:haad} and the discrete gradient method~\eqref{eq:disenen} becomes the Hamiltonian splitting method and the
discrete gradient method for model~\eqref{eq:ffvp}.

\section{Electromagnetic hybrid model}\label{sec:extension}
Here we extend of the aforementioned structure-preserving methods to an electromagnetic hybrid model with Boltzmann electrons and space charge effects proposed in~\cite{vu1996adiabatic}. This model is derived using a temporal WKB approximation when there is a laser with high frequency $w_0$ injected into the plasma, such that numerical simulations on the time scale of the ions can be conducted. 
We assume the vector potential can be written as 
$$
{\mathbf A}({\mathbf x}, t) = \frac{1}{2}\left({\mathbf a}({\mathbf x}, t)e^{-iw_0t} + {\mathbf a}^*({\mathbf x}, t)e^{iw_0t} \right),
$$
where ${\mathbf a} = (a_1, a_2, a_3)^\top = {\mathbf a}_r + \mathrm{i}{\mathbf a}_i$ is complex-valued and is assumed to vary on a time scale much longer than $2\pi/w_0$, and $*$ denotes the conjugate of the complex number.
More details of the derivation can be found in~\cite{vu1996adiabatic}.  After the following normalization, 
\begin{equation*}
\frac{t}{\tilde{t}} = \omega_i^{-1}, \ \frac{\mathbf x}{\tilde{\mathbf x}} = c \omega_i^{-1}, \ \frac{\mathbf v}{\tilde{\mathbf v}} = c, \ \frac{f}{\tilde{f}} = \frac{n_c}{c^3}, \ {\frac{{\mathbf a}}{\tilde{\mathbf a}}} = \frac{cm_i}{e}, \ \frac{\phi}{\tilde{\phi}} = \frac{c^2m_i}{e}, \ \frac{T_e}{\tilde{T_e}} = m_i c^2, \ \omega_i = \sqrt{\frac{n_ce^2}{m_i \epsilon_0}},
\end{equation*}
where $c$ is the speed of light, and $n_c$ is the characteristic density, 
we have the normalized hybrid model~\cite{vu1996adiabatic}
\begin{equation}
\label{eq:electro}
\begin{aligned}
& \frac{\partial f}{\partial t} + {\mathbf v} \cdot \frac{\partial f}{\partial {\mathbf x}} +\left(-Z  \nabla  \phi - \frac{Z^2}{4} \nabla ({\mathbf a} \cdot {\mathbf a}^*) \right)  \cdot \frac{\partial f}{\partial {\mathbf v}} = 0,\\
&  i \epsilon \frac{\partial {\mathbf a}}{\partial t} = - \frac{\epsilon^2}{2} \Delta   {\mathbf a} -  \frac{1}{2}\left(1 -  \epsilon^2 Z^2 \int f \mathrm{d}{\mathbf v} - \epsilon^2 n_e \frac{m_i}{m_e}\right) {\mathbf a}, \quad \omega_i^2 = \frac{e^2 n_c}{m_i\epsilon_0}, \quad \epsilon = \frac{\omega_i}{w_0}\\
& - \Delta \phi =  Z\int f\, \mathrm{d}{\mathbf v} - n_e,
\end{aligned}
\end{equation}
where $Z$ is the ion charge number, $\epsilon$ is very small due to high frequency $w_0$ of the pump wave, and $n_e$ is determined by the potential $\phi$ and ${\mathbf a}$ via the following relations with the given functions $n_0$ and $C$,
$n_e = n_0 e^{\frac{\phi - \frac{m_i}{4m_e}{\mathbf a}\cdot{\mathbf a}^*}{T_e}}$ (isothermal electron case), $n_e = \left( \frac{\phi - \frac{m_i}{4m_e}{\mathbf a}\cdot {\mathbf a}^*}{T_e} \frac{\gamma-1}{\gamma} - C\right)^{\frac{1}{\gamma-1}}, \gamma \neq 1$, (adiabatic electron case).

The equation satisfied by ${\mathbf a}$ is a Schr\"odinger type equation in the form similar to the semiclassical regime~\cite{bao2002time}. Although the small parameter $\epsilon$ in the Schr\"odinger equation introduces oscillations in time and space, the time step size larger than $\epsilon$ is used in~\cite{vu1996adiabatic}. The commonly used numerical scheme for the Schr\"odinger equation in the semiclassical regime is the time splitting spectral method~\cite{bao2002time}, which has the advantage of using large time step and mesh sizes, especially for the computation about 
the observables, such as the term ${\mathbf a}\cdot{\mathbf a}^*$ in this hybrid model.

Regarding the geometric structure,
we propose the following Poisson bracket, which is the sum of the Poisson brackets of the HBS model~\eqref{eq:avp} and the Schr\"odinger equation~\cite{marsden2013introduction} (scaled by $\epsilon$),
\begin{equation}
\label{eq:wkbbracket}
\{\mathcal{F}, \mathcal{G}\}(f,  {\mathbf a}^r, {\mathbf a}^i)  = \int f \left[\frac{\delta \mathcal{F}}{\delta f}, \frac{\delta \mathcal{G}}{\delta f} \right]_{xv} \mathrm{d}{\mathbf x} \mathrm{d}{\mathbf v} +  {\epsilon}  \  \int \frac{\delta F}{\delta {\mathbf a}^r} \cdot \frac{\delta G}{\delta {\mathbf a}^i} - \frac{\delta F}{\delta {\mathbf a}^i}\cdot \frac{\delta G}{\delta {\mathbf a}^r}   \mathrm{d}{\mathbf x},
\end{equation}

The above model~\eqref{eq:electro} can be derived with the above Poisson bracket~\eqref{eq:wkbbracket} and the following Hamiltonian for the isothermal and adiabatic electron cases, respectively
\begin{equation}
\begin{aligned}
\mathcal{H} & = \int  \frac{|{\mathbf v}|^2}{2} f\, \mathrm{d}{\mathbf v}\mathrm{d}{\mathbf x}    +  \int  \frac{|Z{\mathbf a}|^2}{4} f\, \mathrm{d}{\mathbf v}\mathrm{d}{\mathbf x}   +  \frac{1}{4} \int |\nabla a_1|^2 + |\nabla a_2|^2 + |\nabla a_3|^2 \mathrm{d}{\mathbf x}\\
& -  \int \frac{|{\mathbf a}|^2}{4\epsilon}  \mathrm{d}{\mathbf x} - \int T_e n_0 e^{\frac{\phi - \frac{m_i}{4m_e}{\mathbf a}\cdot{\mathbf a}^*}{T_e}} \mathrm{d}{\mathbf x} -  \int \frac{|\nabla \phi|^2}{2} \mathrm{d}{\mathbf x} + \int Zf \phi \,\mathrm{d}{\mathbf x}\mathrm{d}{\mathbf v},\\
\mathcal{H} & = \int  \frac{|{\mathbf v}|^2}{2} f \,\mathrm{d}{\mathbf v}\mathrm{d}{\mathbf x} + \int  \frac{|Z{\mathbf a}|^2}{4} f \,\mathrm{d}{\mathbf v}\mathrm{d}{\mathbf x}   +  \frac{1}{4} \int |\nabla a_1|^2 + |\nabla a_2|^2 + |\nabla a_3|^2 \mathrm{d}{\mathbf x}- \int \frac{|{\mathbf a}|^2}{4\epsilon}  \mathrm{d}{\mathbf x} 
\\
& - \int T_e \left( \frac{\phi - \frac{m_i}{4m_e}{\mathbf a}\cdot {\mathbf a}^*}{T_e} \frac{\gamma-1}{\gamma} - C\right)^{\frac{\gamma}{\gamma-1}} \mathrm{d}{\mathbf x} -  \int \frac{|\nabla \phi|^2}{2} \mathrm{d}{\mathbf x} + \int Zf \phi\, \mathrm{d}{\mathbf x}\mathrm{d}{\mathbf v}.
\end{aligned} 
\end{equation}
The phase-space discretization can be conducted as above through the discretization of the Poisson bracket as~\cite{Qincanonical, GEMPIC}.
The Hamiltonian splitting method~\cite{crouseilles2015hamiltonian, qin2015comment, he2015hamiltonian} gives three explicitly solvable subsystems (or in Fourier space), further details are presented in appendix~\ref{sec:appb}.

\section{Numerical experiments}\label{sec:numsec}

In this section, three numerical experiments: finite grid instability (of an equilibrium), Landau damping (of damping waves), and resonantly excited nolinear ion waves (with non-zero $\phi_0$), are conducted using the normalization~\eqref{eq:llll} to illustrate the conservation properties of the schemes~\eqref{eq:haad}-\eqref{eq:disenen} of the HBS model~\eqref{eq:avp}. The reference density $n_0$ is set to 1, and the unit charge number $Z = 1$. The degree of the shape function is 2, the tolerance for the fixed point iteration is $10^{-12}$, and periodic boundary conditions are used.

\subsection{Finite grid instability}
Finite grid instability in the context of hybrid simulations was firstly reported in~\cite{rambo1995finite} for the hybrid model with quasi-neutrality and Boltzmann electrons. This numerical phenomenon typically arises in standard particle-in-cell methods when the temperature ratio $T_e/T_i \gg 1$, and ions are heated until the ion thermal speed becomes comparable to the ion acoustic speed (and therefore, $T_e/T_i \approx 1$). In ~\cite{rambo1997numerical}, it was noted that finite grid instability also occurs when using traditional particle-in-cell methods for the HBS model, although it is weaker than the hybrid model with quasi-neutrality and Boltzmann electrons. In~\cite{stanier2019fully, fB, fA}, 
the finite grid instability is reduced numerically by using the conservative or bracket-based particle-in-cell methods for the hybrid model with quasi-neutrality and massless electrons.
The finite grid instability of particle-in-cell methods has also been studied in~\cite{huang2016finite, xiao2019structure}, which reveals that the aliased spatial modes are the major cause of the finite grid instability in the particle in cell methods, and geometric particle in cell methods are able to suppress the finite grid instability. 

\begin{figure}
\center{\includegraphics[scale=0.32]{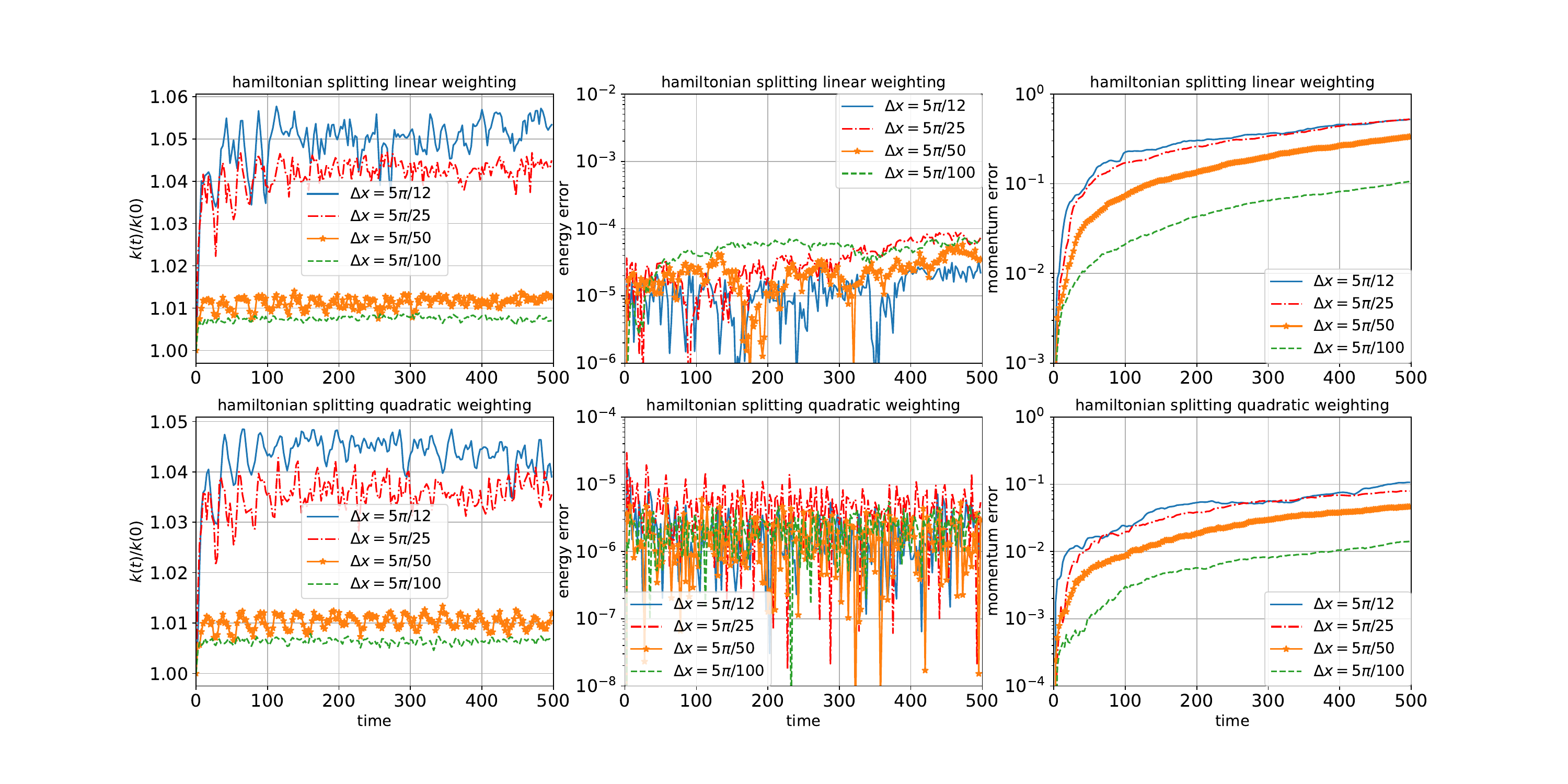}}
\caption{{\bf Finite grid instability of the HBS model by Hamiltonian splitting method.} Time evolution of $k(t)/k(0)$ with $k$ denoting the ion kinetic energy, relative energy error, and momentum error. }
\label{fig:finite}
\end{figure}

\begin{figure}
\center{\includegraphics[scale=0.32]{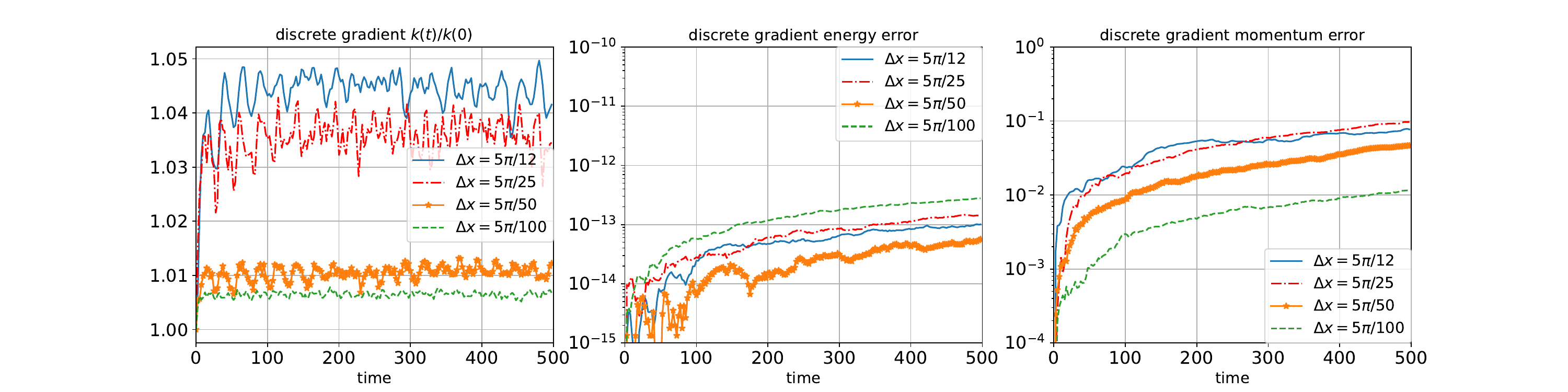}}
\caption{{\bf Finite grid instability of the HBS model by discrete gradient method with quadratic weighting.} Time evolution of $k(t)/k(0)$ with $k$ denoting the ion kinetic energy, relative energy error, and momentum error. }
\label{fig:finitedis}
\end{figure}

In this test, we investigate the finite grid instability using the following initial condition (an equilibrium of~\eqref{eq:avp}) by the numerical simulations conducted with schemes~\eqref{eq:haad} and~\eqref{eq:disenen},
\begin{equation*}
\begin{aligned}
f =  \frac{n_i}{\pi^{\frac{1}{2}} v_T^{\frac{1}{2}}}  \exp\left({\frac{|v - 0.1|^2}{v_T^2} }\right), \quad  T_e = 0.08, \quad v_T = 0.1, \quad n_i = 1.
\end{aligned}
\end{equation*}
Computational parameters include: domain $[0, 5\pi]$, time step size $\Delta t = 0.05$, and particle number per cell $100$. 
We run the simulations with the numerical methods~\eqref{eq:haad} and~\eqref{eq:disenen} with different cell sizes, i.e., $\Delta x = 5\pi/12, 5\pi/25, 5\pi/50, 5\pi/100$, and the results are presented in Fig.~\ref{fig:finite}-\ref{fig:finitedis}. We can see $k(t)/k(0)$ ($k(t) = \frac{1}{2} \sum_{k=1}^{N_p} w_k v_k^2$ the ion kinetic energy) oscillates with time without exhibiting rapid linear growth, as observed in~Fig.~3a of~\cite{rambo1997numerical}. This indicates that 
 the finite grid instability is reduced numerically. 
 As the cell size decreases and the electron Debye length is resolved with higher resolution, the oscillating level of $k(t)/k(0)$ becomes closer to 1, and the momentum error also becomes smaller. The momentum error also depends on the particle number. When there are 100 cells and 2000 particles in each, Hamiltonian splitting method with quadratic weighting gives the momentum error at the level of $10^{-4}$. 
 
As the derivatives of B-splines appear in the schemes~\eqref{eq:haad} and~\eqref{eq:disenen}, second order at least B-spline shape functions should be used.
As the Hamiltonian splitting method~\eqref{eq:haad} (Strang splitting) is symplectic, it has superior long time numerical behaviour, although energy is not conserved exactly (with an relative error about $10^{-5}$) with quadratic weighting. 
In Fig.~\ref{fig:finite}, the relative energy error is also very small ($10^{-4}$) even when the first order B-spline (the derivative of the B-spline is taken as the right derivative). Relative energy error of the discrete gradient method with quadratic weighting is about $10^{-13}$. 

Here we discuss the time step size of the Hamiltonian splitting methods. 
We consider the case with 100 cells in space, time interval $[0,500]$, and quadratic weighting. 
To achieve satisfying results, when $n_i =1$ and particle per cell is 100, the maximum time step size is around 0.5 with energy error at the level of $10^{-4}$. For higher initial densities, such as $n_i = 4, 16$ with  $400, 1600$ particles per cell, the maximum time step sizes giving satisfying results are around $0.3,0.3$ with energy errors at the level of $10^{-4}, 10^{-4}$, respectively.  

\subsection{Landau damping}
\begin{figure}
\center{\includegraphics[scale=0.35]{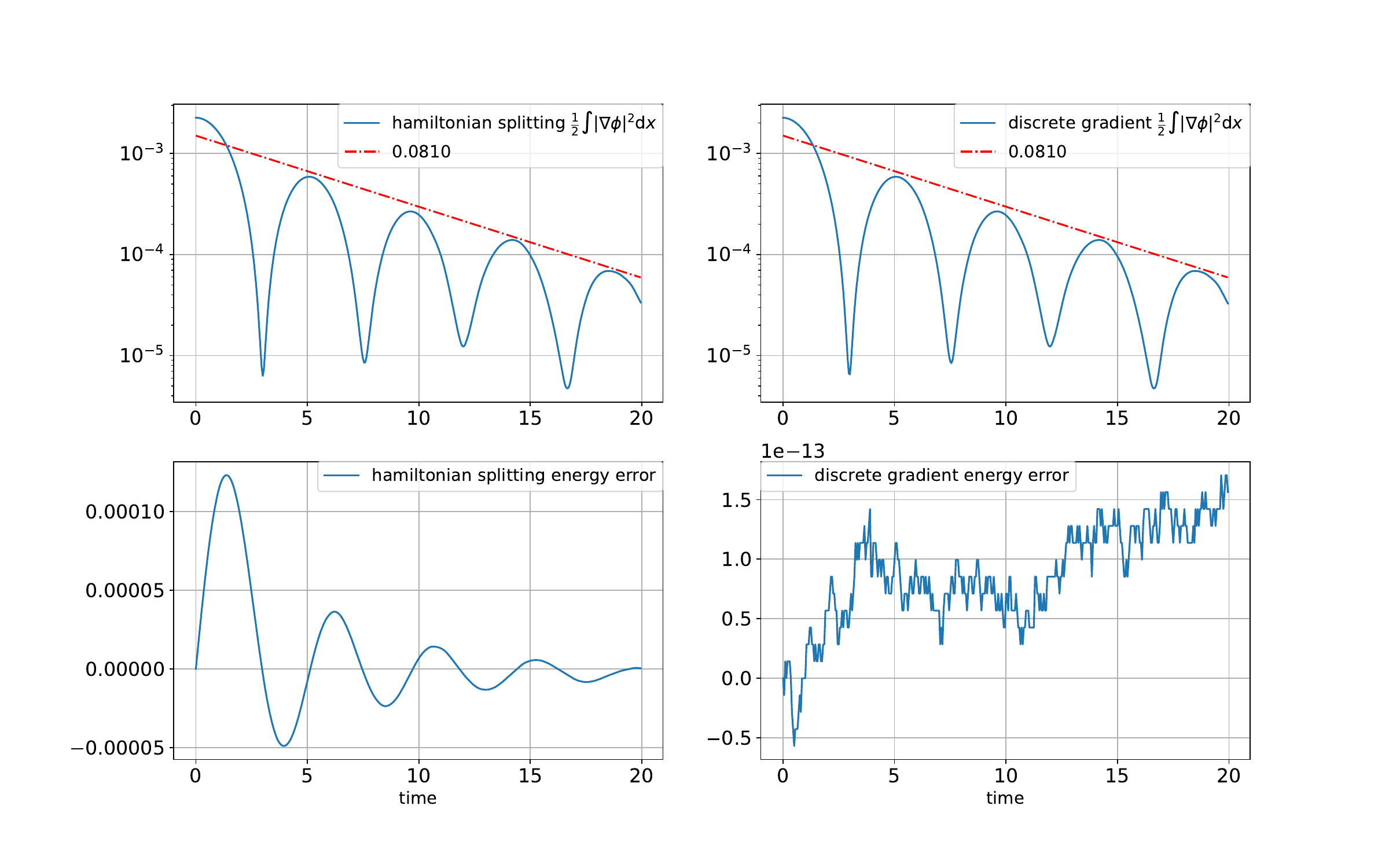}}
\caption{{\bf Linear Landau damping of the HBS model with $T_e=5$ by Hamiltonian splitting and discrete gradient methods.} Time evolution of the electric energy $\frac{1}{2}\int |\nabla \phi|^2  \mathrm{d}x$ and total energy error.}
\label{fig:landaulinear}
\end{figure}

Firstly we simulate the linear ion Landau damping by one-dimensional simulations. The initial distribution function is
\begin{equation*}
\begin{aligned}
f =  \frac{n_i}{\pi^{\frac{1}{2}} v_T^{\frac{1}{2}}} (1 + 0.02\cos(0.25x)) \exp\left({-\frac{v^2}{v_T^2}  }\right). 
\end{aligned}
\end{equation*}
The computational parameters are as follows:
grid number 64, domain size $ [0, 8\pi]$, time step size $\Delta t = 0.05$, final computation time $20$, $v_T = 1.4142$, $n_i = 1$,  and total particle number $10^7$. In this test the quadratic weighting is used.  
See the numerical results with $T_e=5$ in Fig.~\ref{fig:landaulinear} by Hamiltonian splitting method~\eqref{eq:haad} and discrete gradient method~\eqref{eq:disenen}. 
Solving the dispersion relation mentioned in section~\ref{sec:neutrallimit},
$
1 + {k^2 T_e} = \frac{T_e}{T_i}   \mathcal{Z}'\left(\frac{\omega}{kv_T} \right),
$
we find $\omega = 0.6986 - 0.0810\text{i}$ when $k = 0.25$. 
Methods~~\eqref{eq:haad}-\eqref{eq:disenen} give accurate damping rate of the electric energy $\frac{1}{2}\int |\nabla \phi|^2  \mathrm{d}x$. The dispersion relation $
1 = \frac{T_e}{T_i}  \mathcal{Z}'\left(\frac{\omega}{kv_T} \right)
$ of the model with quasi-neutrality and Boltzmann electrons~\cite{drift} in section~\ref{sec:neutrallimit} yields $\omega = 0.7528 - 0.05806\text{i}$, i.e., a slower damping rate.
The energy errors of the schemes~\eqref{eq:haad}-\eqref{eq:disenen} are around $10^{-4}$ and $10^{-13}$, respectively.

\begin{figure}
\center{\includegraphics[scale=0.35]{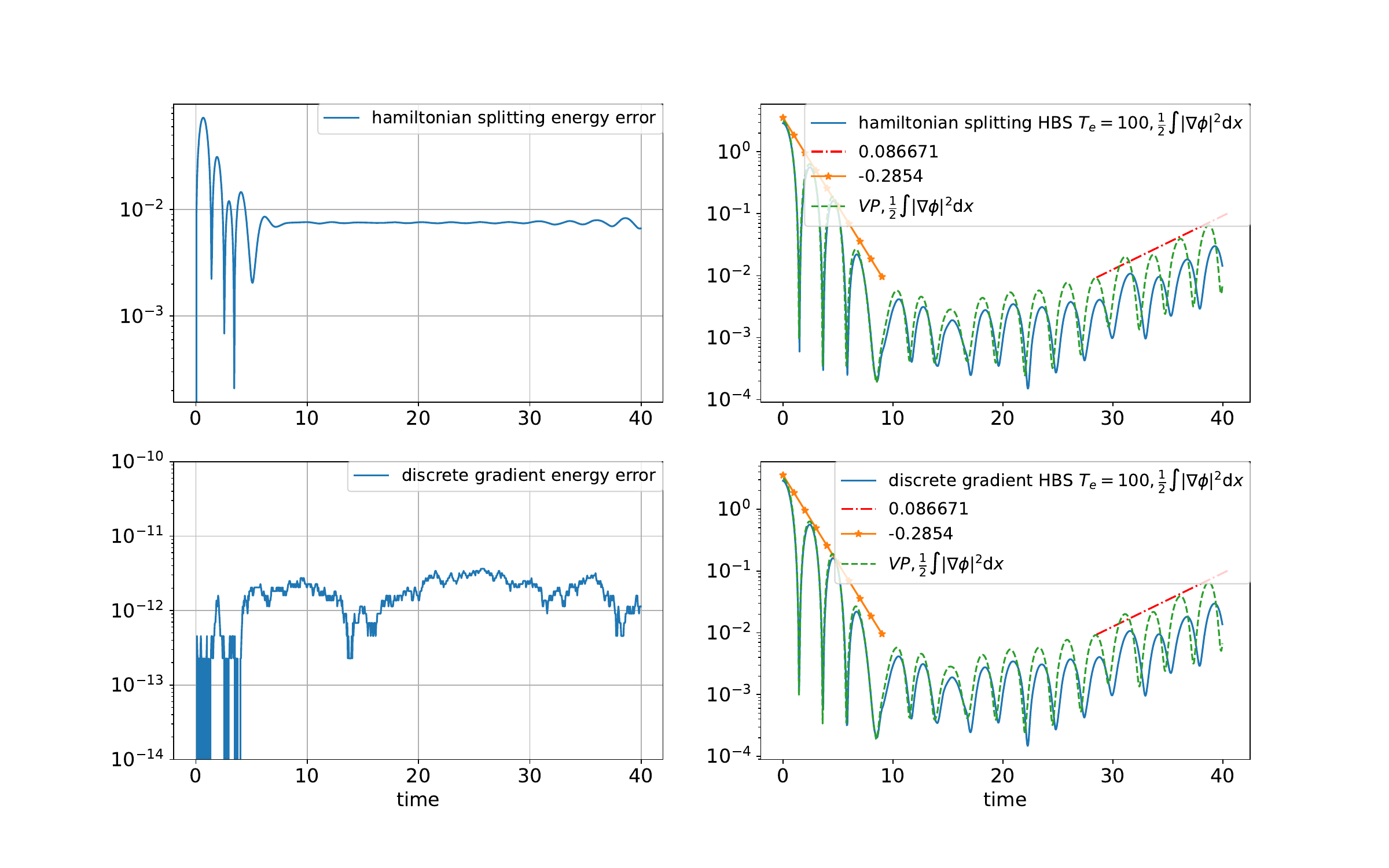}}
\caption{{\bf Nonlinear Landau damping of the HBS model with $T_e=100$ by Hamiltonian splitting and discrete gradient methods.} Time evolution of total energy error and electric energy $\frac{1}{2}\int |\nabla \phi|^2  \mathrm{d}x$. }
\label{fig:nonlinearlandau}
\end{figure}

Then we simulate nonlinear ion Landau damping. The initial distribution function is
\begin{equation*}
\begin{aligned}
f =  \frac{n_i}{\pi^{\frac{1}{2}} v_T^{\frac{1}{2}}} (1 + 0.5\cos(0.5x)) \exp\left({-\frac{v^2}{v_T^2}  }\right). 
\end{aligned}
\end{equation*}
The computational parameters are as follows:
grid number 65, domain size $ [0, 4\pi]$, time step size $\Delta t = 0.05$, final computation time $40$, $v_T = 1.4142, n_i = 1$, and total particle number $10^5$. In this test the quadratic weighting is used.  
See the numerical results with large $T_e=100$ in Fig.~\ref{fig:nonlinearlandau} by Hamiltonian splitting method~\eqref{eq:haad} and discrete gradient method~\eqref{eq:disenen}. Due to the large $T_e$, the term $\exp(\phi/T_e)$ approximates 1, 
makeing the solution of HBS model~\eqref{eq:avp} approximates the solution of the Vlasov--Poisson system (with static electron density as 1). In Fig.~\ref{fig:nonlinearlandau}, we observe the nonlinear Landau damping. The time evolution of energy component $\frac{1}{2}\int |\nabla \phi|^2 \mathrm{d}{x}$ decays exponentially initially, and the decay rate is very close to the decay rate 0.2854 of the Vlasov--Poisson system~\cite{GEMPIC} before time $T=10$. $\frac{1}{2}\int |\nabla \phi|^2 \mathrm{d}{x}$ oscillates when $t \in [10, 30]$, then grows exponentially with time when $t \in [30, 40]$ with a rate close to 0.086671~\cite{GEMPIC} of the Vlasov--Poisson system. For Hamiltonian splitting method, energy error is about $10^{-2}$.
Discrete gradient method gives a smaller total energy error about $10^{-12}$, and similar behavior of electric energy is presented.  The errors of neutrality given by both numerical methods are at the level of iteration tolerance. For the time step size of the Hamiltonian splitting methods, we consider the case with 65 cells in space, time interval $[0,40]$, $10^5$ particles, and quadratic weighting. As $T_e=100$ is large, $\exp(\phi/T_e)$ is close to 1, the numerical stability property is close to the result of~\cite{kormann2021energy}, i.e., the stability condition is around $\Delta t \omega_i < 2$. In order to get the acceptable accuracy, the time step size is usually chosen smaller than 2. When $n_i = 1$, the maximum time step size yielding good numerical behaviour is $\Delta t = 0.4$, resulting in an energy error around $0.45$ after saturation; For $n_i = 4$, $\Delta t = 0.2$  gives an energy error around $0.4$ after saturation, and for $n_i = 16$, $\Delta t = 0.12$  results in an energy error around $0.3$ after saturation. We also consider the case with a small electron temperature $T_e = 1$. In this case, when $n_i = 1$, $\Delta t = 0.5$ gives an energy error around 0.05 after saturation; when $n_i = 4$, $\Delta t = 0.2$ gives an energy error around 0.002 after saturation; when $n_i = 16$, $\Delta t = 0.1$  gives an energy error around 0.004.

\subsection{Simulations with the ponderomotive driving term}\label{subsec:ponderomotive}

\begin{figure}
\center{\includegraphics[scale=0.35]{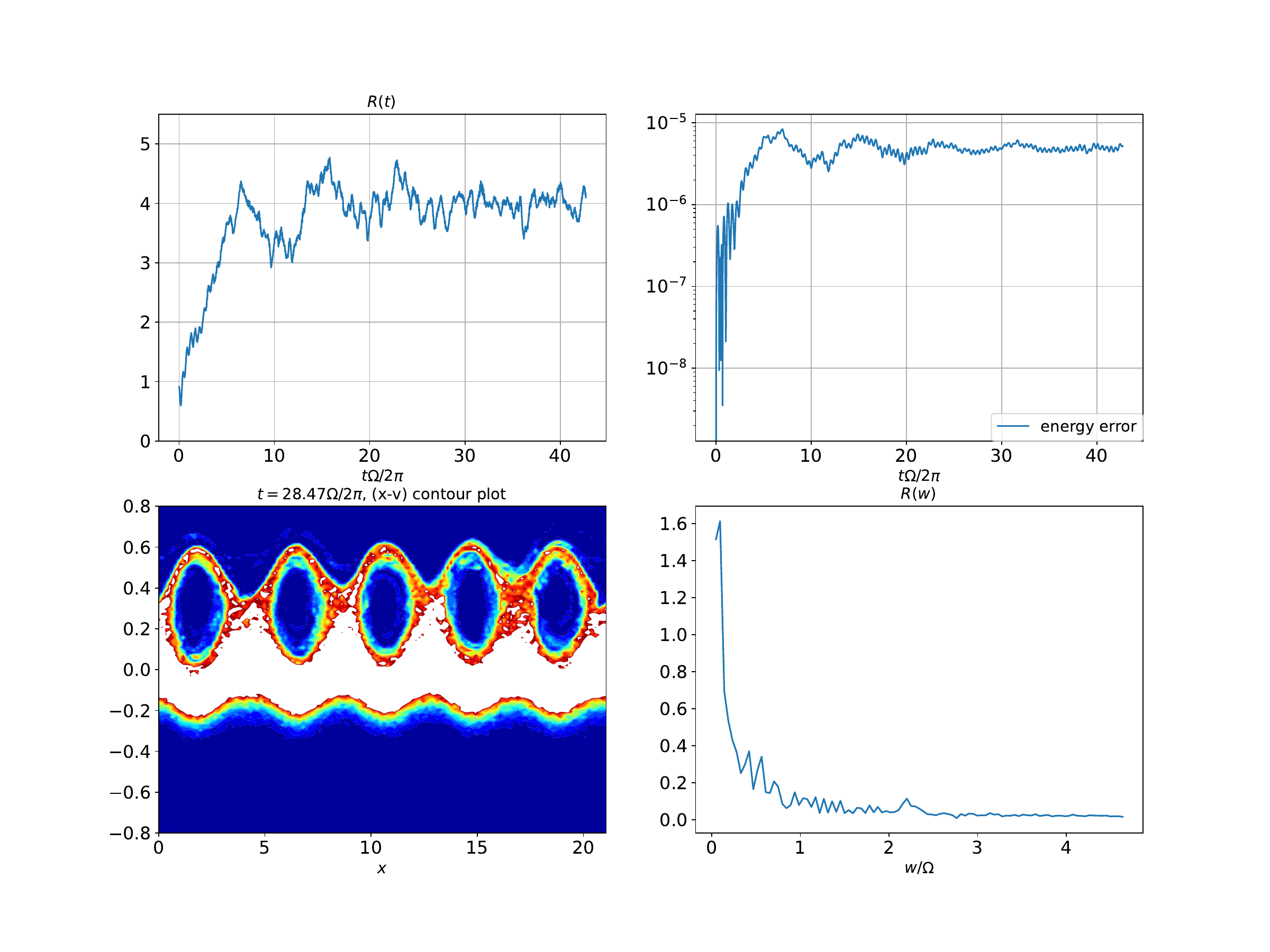}
}
\caption{{\bf Simulations with the ponderomotive driving term by Hamiltonian splitting method.} Time evolutions of $R(t) = \text{max} \frac{\phi}{\tilde{\phi}_0}$ and energy error, the contour plot of the distribution function at time $t=400$, and the fast Fourier transformation of $R(t)$.}
\label{fig:driveha}
\end{figure}

\begin{figure}
\center{\includegraphics[scale=0.35]{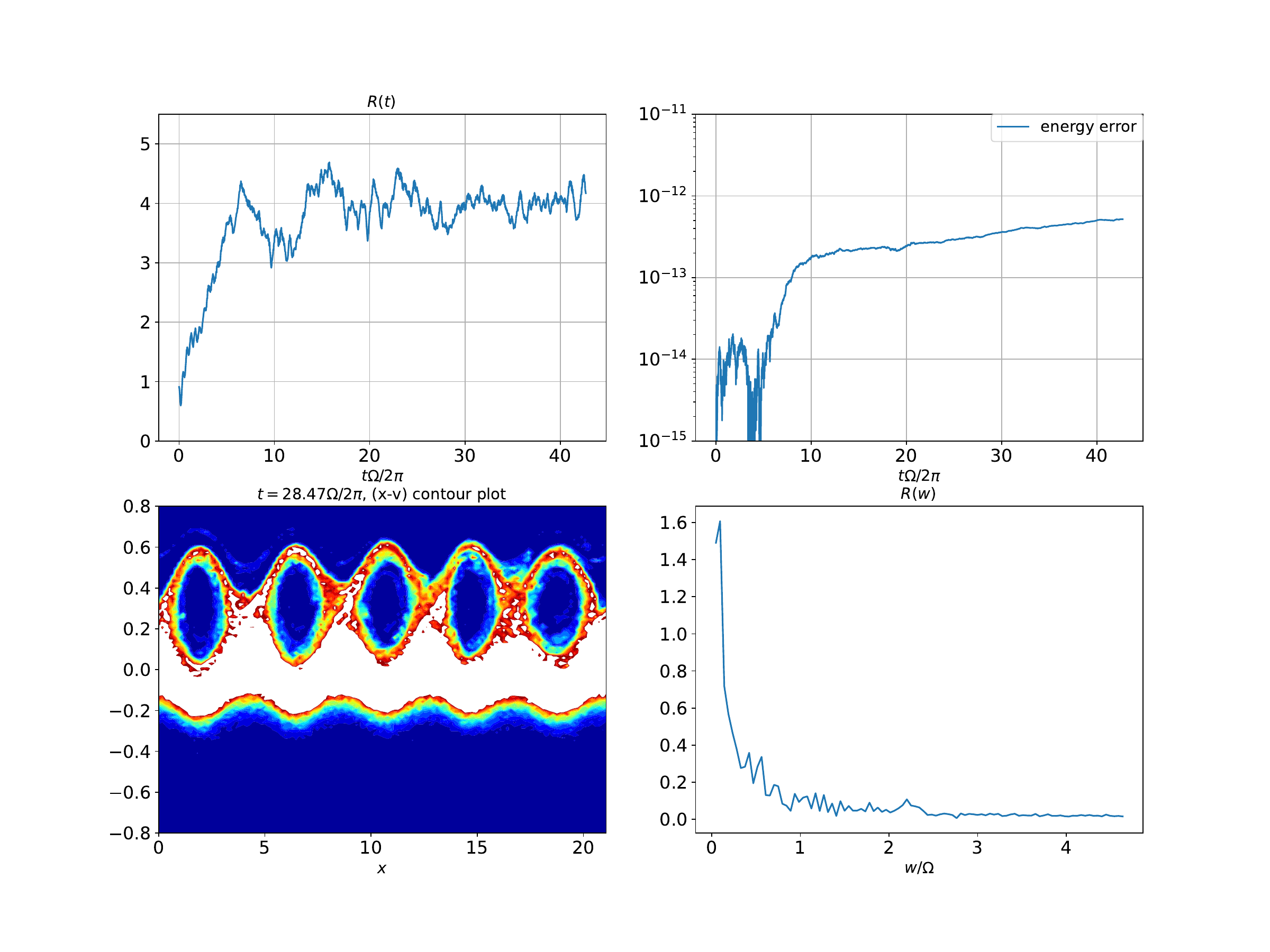}
}
\caption{{\bf Simulations with the ponderomotive driving term by discrete gradient method.} Time evolutions of $R(t) = \text{max} \frac{\phi}{\tilde{\phi}_0}$ and energy error, the contour plot of the distribution function at time $t=400$, and the fast Fourier transformation of $R(t)$.}
\label{fig:drivedis}
\end{figure}

In~\cite{cohen1997resonantly}, the HBS model with a  ponderomotive driving term was numerically solved  to study nonlinear ion acoustic waves. 
Here, following~\cite{cohen1997resonantly}, we conduct a simulation with a non-zero given time dependent function $\phi_0$ in~\eqref{eq:avp}. Specifically, the initial condition and $\phi_0$ are given by
\begin{equation*}
\begin{aligned}
f =  \frac{n_i}{\pi^{\frac{1}{2}} v_T^{\frac{1}{2}}} \exp\left({-\frac{v^2}{v_T^2}  }\right), \quad \phi_0 = \tilde{\phi}_0 \cos(\Omega t - k x),
\end{aligned}
\end{equation*}
where $n_i = 1, v_T = \frac{\sqrt{2}}{10}$, $\tilde{\phi}_0 = 0.05T_e$, $\Omega = 0.4472$, $k = 1.49$. Other computational parameters are: grid number 64, domain size $ [0, \frac{10\pi}{k}]$, time step size $\Delta t = 0.1$, final computation time $600$,  $T_e = 0.1125$, and total particle number $10^6$. In this test the quadratic weighting is used. Since $\phi_0$ is time dependent, the Hamiltonian system~\eqref{eq:spdeha} is a non-autonomous Hamiltonian system, for which we use the the technique of extending the dimension~\cite{zhou2017explicit}. From Fig.~\ref{fig:driveha}-\ref{fig:drivedis}, we can see that the peak value of the response function $R(t) = \text{max}_{x}\frac{\phi}{\phi_0}$ is around 5, which is consistent
with the result in~\cite{cohen1997resonantly}. 
There is a rapid oscillation at $2.1\Omega$ and a slow modulation at $0.04\Omega$ in the fourth figure obtained by the fast Fourier transformation of $\text{max}_x\frac{\phi}{\phi_0}$ in time, which are close to the results in~\cite{cohen1997resonantly} with $1.85 \Omega$ (fast) and $0.15\Omega$ (slow).
The Hamiltonian splitting method and discrete gradient method give the energy errors around $10^{-5}$ and $10^{-12}$, respectively. And when $t=400$, both methods give similar 5 vortices in phase-space contour plot, due to the the driving force is the fifth mode, i.e., $k = 5 \frac{2\pi}{L}$, where $L$ is the domain size. 

Regarding the time step size for the Hamiltonian splitting methods, we use the same parameters as above but vary the initial density and time step size. When the initial density $n_i =  1$, the maximum step size for giving good numerical behaviour is $\Delta t = 1.5$, which gives an energy error around 0.003; When the initial density $n_i = 4$, the maximum step size $\Delta t = 1$ gives an energy error around 0.02; When the initial density $n_i = 16$, the maximum step size $\Delta t = 0.5$ gives an energy error around 0.008.

\section{Conclusion}\label{sec:conclusion}
In this paper, we explore the structure-preserving discretizations of the electrostatic hybrid  plasma model with Boltzmann electrons and space-charge effects. These discretizations are derived by discretizing either the variational action integral or the Poisson bracket combined with the Hamiltonian splitting methods~\cite{crouseilles2015hamiltonian, qin2015comment, he2015hamiltonian} in time. 
Discrete gradient methods~\cite{mclachlan1999geometric} are employed to conserve energy exactly. The geometric structure and numerical discretization of the electromagnetic hybrid model~\cite{vu1996adiabatic} are detailed in Section~\ref{sec:extension}.

For discretizing the field functions, 
the finite element methods~\cite{GEMPIC} or Fourier spectral methods~\cite{pinto2021geometric} can be used, while the distribution function can be discretized by the delta functions. Additional details can be found in the appendix~\eqref{sec:appa}. The cases with other kinds of boundary conditions and further exploration (such as the physical application and the time and mesh size strategy) of the electromagnetic hybrid model can be considered in future works.

%
%

\section*{Acknowledgements}
The simulations in this work were performed on Max Planck Computing \& Data Facility (MPCDF). The author would like to thank
anonymous reviewers for many helpful comments for improving this paper.  Special thanks to B. I. Cohen for the help of the simulation parameters used in~\ref{subsec:ponderomotive}. The author would like to acknowledge P. J. Morrison, S. Possanner, and E. Sonnendr\"ucker for their kind discussions of this work.

\section*{Declaration of interests}
The authors report no conflict of interest.

\appendix
\section{Discretization with finite element method}\label{sec:appa}

Distribution function $f$ is approximated using $\delta$ functions, i.e, 
$$
f({x}, {v}, t) \approx f_h({x}, {v}, t) = \sum_{k=1}^{N_p} w_k \delta ({x} - {x}_k) \delta ({v} - {v}_k),
$$
where $N_p$ is the total particle number, and $w_k$, ${x}_k$, and ${v}_k$ are the weight, position, and velocity for $k$-th particle. 
We discretize $\phi$ by finite element method, i.e. 
$$
\phi_h = {\boldsymbol \Lambda} \cdot {\boldsymbol \phi},
$$
where the vectors  ${\boldsymbol \Lambda}$ and ${\boldsymbol \phi}$ contain all basis functions and finite element coefficients.  The Poisson--Boltzmann equation is discretized  in weak formulation as 
\begin{equation*}
\underbrace{\int \partial_x \phi_h   \partial_x \Lambda_i \mathrm{d}{x}}_{= \mathbb{M} {\boldsymbol \phi}} + \underbrace{\int n_0 \exp\left({\frac{\phi_h - \phi_{0,h}}{T_e}}\right) \Lambda_i \mathrm{d}{x}}_{\approx \sum_{j=1}^N w_j n_0({ x}_j) \exp\left({\frac{\phi_h( {x}_j) - \phi_{0,h}( {x}_j)}{T_e({x}_j)}}\right) \Lambda_i({x}_j)} = Z \sum_{k=1}^{N_p} w_k \Lambda_i({x}_k),
\end{equation*}
where ${x}_j$ is the $j$-th quadrature point, $w_j$ is the corresponding quadrature weight, and $\mathbb{M}_{ij} = \int  \partial_x \Lambda_i  \partial_x \Lambda_j \mathrm{d}{x}$.\\
We approximate variational action integral~\eqref{eq:vpai} as 
\begin{equation}\label{eq:femva}
\mathcal{A}_h = \sum_{k=1}^{N_p} w_k \left( \frac{\dot{x}_k^2}{2} - Z \phi_n({x}_k)\right) + {\boldsymbol \phi}^\top \mathbb{M} {\boldsymbol \phi} + \sum_{j=1}^N w_j n_0({x}_j) T_e({x}_j) \exp\left(\frac{\phi_h({x}_j) - \phi_{0,h}(x_j)}{T_e({x}_j)}\right).
\end{equation}
Hamiltonian is discretized as 
\begin{equation}\label{eq:femap}
H = \sum_{k=1}^{N_p} w_k \frac{{v}_k^2}{2} + Z w_k {\boldsymbol \Lambda}({x}_k) \cdot {\boldsymbol \phi} - \sum_{j=1}^N w_j T_e({x}_j) n_0({x}_j) \exp\left({\frac{\phi_h( {x}_j) -\phi_{0,h}( {x}_j) }{T_e({x}_j)}}\right)  -  \frac{{\boldsymbol \phi}^\top \mathbb{M} {\boldsymbol \phi}}{2}.
\end{equation}
As~\cite{Qincanonical, GEMPIC}, the bracket is discretized as 
\begin{equation}\label{eq:femap}
\{F, G\}_h = \sum_{k=1}^{N_p} \frac{1}{w_k} \left(\partial_{x_k}F    \partial_{v_k}G -  \partial_{x_k}G    \partial_{v_k}F \right).
\end{equation}
Both the variations of \eqref{eq:femva} and the discrete Poisson bracket~\eqref{eq:femap} with Hamiltonian~\eqref{eq:femap} give the following Hamiltonian ODE,
\begin{equation*}
\dot{x}_k = \frac{1}{w_k} \partial_{v_k}H, \quad \dot{v}_k = - \frac{1}{w_k} \partial_{x_k}H,
\end{equation*}
Similarly, for the cases with periodic boundary conditions, we can prove that the neutrality condition holds in a weak sense, i.e., 
$$ \sum_{j=1}^N w_j n_0({ x}_j) \exp\left({\frac{\phi_h( {x}_j) - \phi_{0,h}( {x}_j)}{T_e({x}_j)}}\right) \Lambda_i({x}_j) = Z \sum_{k=1}^{N_p} w_k \Lambda_i({x}_k), \quad \forall i = 1, \cdots, N.$$

\section{Hamiltonian splitting method for the electromagnetic hybrid model~\cite{vu1996adiabatic}}\label{sec:appb}
We take the isothermal electron case with the following energy for an example
\begin{equation}
\begin{aligned}
\mathcal{H} & = \frac{1}{2} \int  |{\mathbf v}|^2 f\, \mathrm{d}{\mathbf v}\mathrm{d}{\mathbf x} + \frac{1}{4} \int  |Z{\mathbf a}|^2 f\, \mathrm{d}{\mathbf v}\mathrm{d}{\mathbf x}   +  \frac{1}{4} \int |\nabla a_1|^2 + |\nabla a_2|^2 + |\nabla a_3|^2 \mathrm{d}{\mathbf x} \\
& - \frac{1}{4\epsilon} \int |{\mathbf a}|^2  \mathrm{d}{\mathbf x} - \int T_e n_0 e^{\frac{\phi - \frac{m_i}{4m_e}{\mathbf a}\cdot{\mathbf a}^*}{T_e}} \mathrm{d}{\mathbf x} - \frac{1}{2} \int |\nabla \phi|^2 \mathrm{d}{\mathbf x} + \int Zf \phi\, \mathrm{d}{\mathbf x}\mathrm{d}{\mathbf v}.
\end{aligned}
\end{equation}
We split the Hamiltonian into the following 3 parts and get the corresponding explicitly solvable subsystems,
\begin{equation*}
\begin{aligned}
\mathcal{H} & = \underbrace{\frac{1}{2} \int  |{\mathbf v}|^2 f\, \mathrm{d}{\mathbf v}\mathrm{d}{\mathbf x} + \frac{1}{4} \int |\nabla a_1|^2 + |\nabla a_2|^2 + |\nabla a_3|^2 \mathrm{d}{\mathbf x} }_{H_{1}}\\
& + \underbrace{\frac{1}{4} \int  |Z{\mathbf a}|^2 f\, \mathrm{d}{\mathbf v}\mathrm{d}{\mathbf x} - \frac{1}{4\epsilon} \int |{\mathbf a}|^2  \mathrm{d}{\mathbf x}  }_{H_{2}}\\
& \underbrace{- \int T_e n_0 e^{\frac{\phi - \frac{m_i}{4m_e}{\mathbf a}\cdot{\mathbf a}^*}{T_e}} \mathrm{d}{\mathbf x} - \frac{1}{2} \int |\nabla \phi|^2 \mathrm{d}{\mathbf x} + \int Zf \phi\, \mathrm{d}{\mathbf x}\mathrm{d}{\mathbf v}}_{H_3}.
\end{aligned} 
\end{equation*}
\noindent{\bf Subsystem $H_{1}$} Corresponding subsystem is 
\begin{equation*}
 \frac{\partial f}{\partial t} + {\mathbf v} \cdot \frac{\partial f}{\partial {\mathbf x}} = 0,\quad  i \epsilon \frac{\partial {\mathbf a}}{\partial t} = - \frac{\epsilon^2}{2} \Delta   {\mathbf a},
\end{equation*}
where the first equation is an explicitly solvable transport equation, and the second equation can be solved explicitly in Fourier space. 

\noindent{\bf Subsystem $H_{2}$} Corresponding subsystem is 
\begin{equation*}
 \frac{\partial f}{\partial t} - \frac{Z^2}{4} \nabla ({\mathbf a} \cdot {\mathbf a}^*)  \cdot \frac{\partial f}{\partial {\mathbf v}} = 0,\quad  i \epsilon \frac{\partial {\mathbf a}}{\partial t} = -  \frac{1}{2} \left(1 - \epsilon^2 Z^2 \int f\, \mathrm{d}{\mathbf v}\right){\mathbf a},
\end{equation*}
where the ${\mathbf a} \cdot {\mathbf a}^*$ is preserved by the second equation and $\int f\, \mathrm{d}{\mathbf v}$ is preserved by the first equation, which make the two equations explicitly solvable.

\noindent{\bf Subsystem $H_{3}$} Corresponding subsystem is 
\begin{equation*}
\frac{\partial f}{\partial t} -Z  \nabla  \phi  \cdot \frac{\partial f}{\partial {\mathbf v}} = 0,\quad  i \epsilon \frac{\partial {\mathbf a}}{\partial t} =  \frac{1}{2} \epsilon^2 n_e \frac{m_i}{m_e}{\mathbf a}, \quad  -\Delta \phi =  Z\int f\, \mathrm{d}{\mathbf v} - n_e,
\end{equation*}
where ${\mathbf a} \cdot {\mathbf a}^*$ is preserved by the second equation, and the ion density $\int f\, \mathrm{d}{\mathbf v}$ is preserved by the first equation. As $n_e$ depends on $\phi$ and ${\mathbf a} \cdot {\mathbf a}^*$, and ${\mathbf a} \cdot {\mathbf a}^*$ and $\int f\, \mathrm{d}{\mathbf v}$ are not changed in this sub-system, $\phi$ and $n_e$ are preserved by this subsystem. We only need to solve the third equation for obtaining $\phi$ once in each time step, and the first and second equations are explicitly solvable.

\bibliographystyle{plain}
\bibliography{hybrid}

\end{document}